\documentclass[12pt]{article}
\usepackage{graphics}
\usepackage{natbib}
\usepackage{epsfig}
\usepackage{amssymb}

\usepackage{color}

\def\qed{\hfill  \framebox(5,5){}}
\def\para{\vspace{2mm}}

\newtheorem{theorem}{Theorem}[section]
\newtheorem{proposition}[theorem]{Proposition}
\newtheorem{corollary}[theorem]{Corollary}
\newtheorem{lemma}[theorem]{Lemma}
\newtheorem{remark}[theorem]{Remark}
\newtheorem{definition}[theorem]{Definition}
\newtheorem{example}[theorem]{Example}

\def\card{{\rm Card}}

\def\cc{{\cal C}}
\def\conch{\mathfrak{C}(\cc,A,d)}
\def\bb{\mathfrak{B}(\cc,A,d)}

\def\K{\Bbb  K}
\def\lmas{{\cal   L}^{+}}
\def\lmenos{{\cal   L}^{-}}
\def\eq{{\rm Eq}}

\begin{document}

\title{
 An Algebraic Analysis of   Conchoids  to Algebraic
Curves.
 \footnote{Both authors supported by the Spanish `` Ministerio de
Educaci\'on y Ciencia" under the Project MTM2005-08690-C02-01}}
\author{
$\begin{array}{ccc}\mbox{J. Rafael
Sendra} & &  \mbox{Juana Sendra} \\
\mbox{Dep. de Matem\'aticas} & &  \mbox{Dep. de Matem\'aticas}\\
\mbox{Universidad de Alcal\'a} & &
\mbox{E.U.I.T. Telecomunicaci\'on} \\
 \mbox{Alcal\'a de Henares, Madrid, Spain} & & \mbox{Univ. Polit\'ecnica de
 Madrid, Spain}\\
 \mbox{{\tt rafael.sendra@uah.es}} & & \mbox{{\tt jsendra@euitt.upm.es}}
  \\
\end{array}$ \and  $\begin{array}{c}
\end{array}$
}
\date{}
\maketitle
\begin{abstract}
We study conchoids to algebraic curve from the perspective of
algebraic geometry, analyzing their main algebraic properties. We
introduce the formal definition of conchoid of an algebraic curve by
means of incidence diagrams.  We prove that, with the exception of a
circle centered at the focus and taking $d$ as its radius, the
conchoid is an algebraic curve having at most two irreducible
components.   In addition, we introduce the notions of special and
simple components of a conchoid.   Moreover we state that, with the
exception of lines passing through the focus, the conchoid always
has at least one simple component and that, for almost every
distance, all the components of the conchoid are simple. We state
that, in the reducible case, simple conchoid components are
birationally equivalent to the initial curve, and we show how
special components can be used to decide whether a given algebraic
curve is the conchoid of another curve.
\end{abstract}

\section{Introduction.}
\begin{figure}[ht]
\begin{center}
\centerline{
\psfig{figure=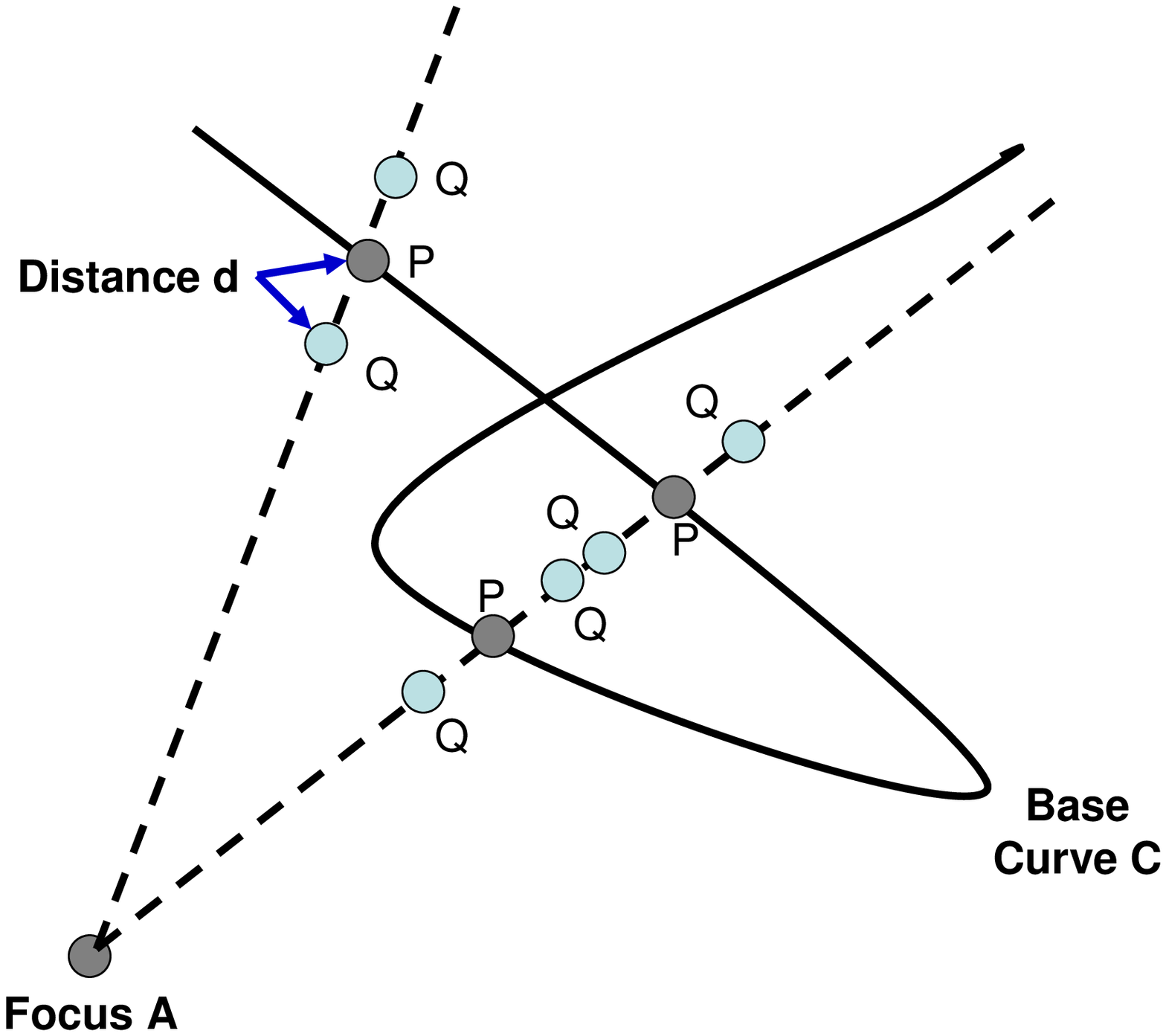,width=7cm,height=8cm,
angle=270}
\psfig{figure=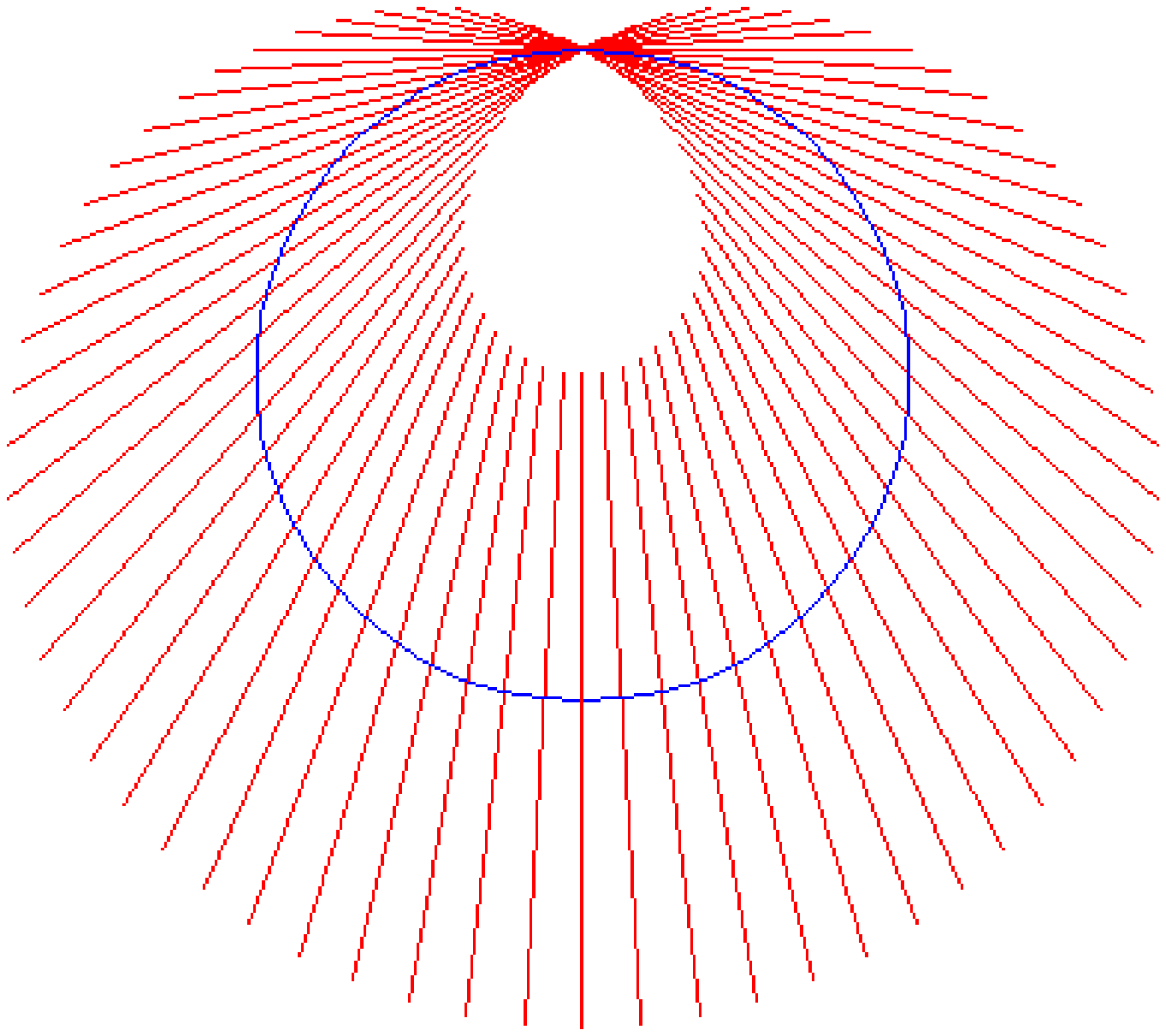,width=5cm,height=5cm,angle=270} }
\end{center}
\caption{{\sf Left:} Conchoid Geometric Construction. \newline
\hspace*{1.5 cm} {\sf Right:} Conchoid of a circle with focus on
it.}
 \label{fig-idea-intuitica}
\end{figure}

A conchoid is a curve derived from a fixed point, another curve, and
a  length in the following way. Let $\cc$ be a plane curve (the base
curve),  $A$ a fixed point in the plane (the focus), and $d$ a
non-zero fixed field element (the distance). Then, we consider the
set of all points in the plane for which there exists a point
$P\in\cc$ such that the distance between $P$ and $Q$ is $d$, being
$A,P,Q$  collinear (see Fig. \ref{fig-idea-intuitica}). Such a
geometric locus will be called the {\it conchoid of $\cal C$ from
the focus $A$ at distance $d$}.

\para

The two classical and most famous conchoids are the {\it Conchoid of
Ni\-co\-me\-des} (see Example \ref{Example-3} and Fig.
\ref{fig-nicomedes}) and the {\it  Lima\c{c}on of Pascal} (see
Example \ref{Example-2}, Fig. \ref{fig-idea-intuitica} right, Fig.
\ref{fig-caracol}) that appear when the base curve $\cc$ is a line
and a circle, respectively.  Conchoid of Nicomedes was introduced by
Nicomedes, around 200 {\sc b.c.}, to solve the problems of doubling
the cube and of trisecting an angle.

\para

Conchoids play an important role in many applications as
construction of buildings (one can already find  specific methods
for producing Lima\c{c}ons in Albert D\"{u}rer's {\it Underweysung
der Messung}), astronomy (see \cite{Kerrick}),  electromagnetic
research (see \cite{Weigan}), physics (see \cite{Szmulowicz96}),
optics (see \cite{Azzam}), engineering in medicine and biology (see
\cite{Menschik97}, \cite{MyungJin04}), mechanical in fluid
processing (see \cite{Sultan}), etc.

\para

Although conchoids have been extensively used and applied in
different areas, a deep theoretical analysis of the concept and its
main properties is missing; at least from our point of view.  In
this paper, we consider conchoids from the perspective of algebraic
geometry, and we study  their main algebraic properties, with the
aim of building a solid bridge from theory to practice that can be
used for further theoretical and applied developments. More
precisely, we introduce the formal definition of conchoid of an
algebraic curve, over an algebraically closed field of
characteristic zero, by means of incidence diagrams. We also
introduce the notion of generic conchoid, and we show how
elimination theory techniques, as Gr\"obner bases, can be applied to
compute conchoids. We prove that, with the exception of a circle
centered at the focus and taking $d$ as its radius, the conchoid is
an algebraic curve having at most two irreducible components. Note
that for the particular circle, mentioned above, the conchoid
consists in two components: a circle of radius $2d$ and the
zero-dimensional set formed by the focus. In addition, we introduce
the notions of special and simple components of a conchoid.
Essentially, a component of a conchoid is special if its points are
generated for more than one point of the original curve. This
phenomenon appears when one computes conchoids of conchoids.
Moreover we state that, with the exception of lines passing through
the focus, the conchoid always has at least one simple component.
Furthermore we prove that, for almost every distance and with the
exception of lines passing through the focus, all the components of
the conchoid are simple. Simple and special components play an
important role in the study of conchoids. On one hand, simple
components are related to the birationality of the maps in the
incidence diagram (for instance, if a conchoid has two components,
its simple components are birationally equivalent to the initial
curve) and, on the other hand, special components can be used to
decide whether a given algebraic curve is the conchoid of another
curve. A similar behavior of simple components holds for offsets
(see \cite{SS99}), and have allowed us to provide formulas for the
genus (see \cite{ASS97}). We plan to investigate this, for the case
of conchoids, in our future research.

\para

The paper is structured as follows. In Section
\ref{sec-def-conchoid}, we formally introduce the notion of conchoid
of an algebraic curve. In Section \ref{sec-basic-prop-conchoids} we
state the basic algebraic properties of conchoids. Section
\ref{sec-simple-special}, we introduce the notion of simple and
special component and we state its main properties. Section
\ref{sec-role} studies how simple components are related to the
birationality of the maps in the incidence diagram, and Section
\ref{sec-detecting} shows how special components can be applied to
detect whether a curve is the conchoid of another curve.

\section{Definition of Conchoid.}\label{sec-def-conchoid}
Let  $\K$ be an algebraically closed field of characteristic zero.
We consider $\K^2$ as the metric affine space induced by the inner
product $B((x_{1},x_{2}),(y_{1},y_{2}))= x_{1}y_{1}+x_{2}y_{2}$. In
this context, the circle of center $(a_{1},a_{2})\in \K^2$ and
radius $d\in \K$ is the plane curve defined by
$(x_{1}-a_{1})^2+(x_{2}-a_{2})^2=d^2$.
We will say that the distance between the points
$\bar{x}$, $\bar{y} \in \K^2$ is $d\in \K^*$ if $\bar{y}$ is on the
circle of center $\bar{x}$ and radius $d$ (notice that the distance
is hence defined up to the sign). On the other hand, if $\bar{x}\in
{\Bbb  K}^{2}$ is not isotropic we denote by $\| \bar{x}\|$ any of
the numbers such that $\|\bar{x}\|^2=B(\bar{x},\bar{x})$, and if
$\bar{x}\in {\Bbb  K}^2$ is isotropic, then $\|\bar{x}\|=0$. In this
paper we usually work with both solutions of
$\|\bar{x}\|^2=B(\bar{x},\bar{x})$. For this reason we use the
notation $\pm \|\bar{x}\|$.

\para

In this situation,  let $\cc$ be the affine irreducible plane curve
 defined by the irreducible polynomial
$f(y_1,y_2)\in \K[y_1,y_2]$, let $d\in \K^*$ be a non-zero field
element, and let $A=(a,b)\in \K^2$. In order to get a formal
definition of the conchoid, one introduces the following incidence
diagram:
\[\begin{array}{c}
 \fbox{
 $ \hspace*{5 mm}
 \begin{array}{ccc} \\
& \bb
 \subset  {\Bbb  K}^{2} \times {\Bbb  K}^{2} \times {\Bbb  K} & \\
& \hspace*{0.5 cm} \mbox{ $\begin{array}{c} \vspace*{-2 mm}
\pi_{1}  \\
\\
\end{array}$ \hspace{-2 mm} {\Huge $\swarrow$}
} \hspace{2 mm} \mbox{ {\Huge $\searrow$} $\begin{array}{c}
\vspace*{-2 mm}
\pi_{2} \\
\\
\end{array}$} &   \\
& \hspace*{-5 mm} \pi_{1}(\bb)
 \subset
{\Bbb  K}^{2} \hspace{1 cm}    {\cal C}
 \subset {\Bbb  K}^{2} & \\ & &  \\
\end{array}$ \newline
\mbox{\sf (Incidence Diagram)}}
\\ \\
\end{array} \]
 where $\bb$ is the algebraic  set of $\K^2\times \K^2\times \K$  defined
 as
\[
\left\{ (\bar{x}, \bar{y},w) \in {{\Bbb  K}}^{2} \times {{\Bbb
K}}^{2}\times {{\Bbb  K}} \left/
\begin{array}{l}
f(y_{1},y_{2})=0   \\
(x_{1}-y_{1})^2+ (x_{2}-y_{2})^2=d^2 \\
 (y_{2}-b)\cdot (x_{1}-y_{1})- (y_{1}-a)\cdot (x_{2}-y_{2})=0 \\
 w\cdot((y_{1}-a)^2+(y_{2}-b)^2)=1
 \end{array}
\right \} \right.
\]
and where  $\pi_{1}$, $\pi_{2}$ are the natural projections
\[
\begin{array}{llllllllll}
\pi_{1}:&  {\Bbb  K}^{2} \times {\Bbb K}^{2}\times {{\Bbb
K}}&\longrightarrow {\Bbb K}^{2},& & & \pi_{2}: & {\Bbb  K}^{2}
\times {\Bbb K}^{2}\times {{\Bbb
K}}&\longrightarrow {\Bbb K}^{2}\\
& (\bar{x}, \bar{y},w) & \longmapsto \bar{x} & & & & (\bar{x},
\bar{y},w) & \longmapsto \bar{y}.
\end{array}
\]
Observe that the first equation in $\bb$ corresponds to $\cc$, the
second and third guarantee that the distance between $\bar{x}$ and
$\bar{y}$ is $d$ and that $\bar{x},\bar{y}, A$ are collinear. The
last equation excludes isotropic points on $\cc$; i.e. point
$\bar{y}\in \cc$ such that $(y_{1}-a)^2+(y_{2}-b)^2=0$. Note that
this  phenomenon occurs, for instance, taking
$f(y_1,y_2)=(y_1-a)+\sqrt{-1} (y_2-b)$. This type of situation will
be analyzed in the next section. Then, we introduce the conchoid as
follows.


\begin{definition}\label{def-conchoide} Let $\cc$ be an affine irreducible plane curve,
$d\in \K^*$, and  $A\in \K^2$. We define the {\sf conchoid of the
base curve
 ${\cal C}$ from the focus
$A$ and distance $d$} as the algebraic Zariski   closure in ${\Bbb
K}^{2}$ of $\pi_1(\bb)$, and we denote it by $\conch$. That is,
\[ \conch= \overline{\pi_{1}(\bb)}. \]
\end{definition}

\begin{remark}\label{remark-a-la-def-conch}
Observe that:
\begin{itemize}
\item[(1)] {\sf [Extension of the definition].} In Definition \ref{def-conchoide} we have considered irreducible
curves. The same reasoning can be done for reducible curves,
introducing the conchoid as the union of the conchoids of the
irreducible components. \item[(2)] {\sf [False components].} Let
$\cc$ be a line (different to $y_1\pm \sqrt{-1} y_2=0$), and $A\in
\cc$ (similarly for any other irreducible curve). Then, following
the ``intuitive" geometric description of the conchoid, one might
claim that $\conch$ consists in $\cc$ and a circle centered at $A$
and radius $d$. However, attending to Definition
\ref{def-conchoide}, the circle (let us call it $\cal D$) is not a
component of the conchoid. Note that, although for every $\bar{x}\in
{\cal D}$ the pair $(\bar{x},A)$ satisfies the three equations of
$\bb$ it does not satisfies the last one. \item[(3)] {\sf
[Computation of the Conchoid].} Let $I$ be ideal in
$\K[\bar{x},\bar{y},w]$ generated by the polynomials defining
$\mathfrak{B}(\cc,A,d)$. Then, by the Closure Theorem (see
\cite{Cox} p. 122), one has that $\conch=V(I\cap \K[\bar{x}])$.
Hence elimination theory techniques, as Gr\"obner bases, provide the
conchoid. \item[(4)] {\sf [Generic Conchoid].} Reasoning as in
Section  2 in \cite{SSS06}, one may introduce the notion of generic
conchoid.  Let us consider  $d$ as a new variable. Now, $\bb$ is
seen as an algebraic set in $\K^2\times \K^2\times \K\times \K$; we
denote it by $\bb_G$. Then, the generic conchoid is defined as
\[ \conch_G= \overline{\pi_{1}(\bb_G)}. \] Now, if $I_G$ is the
ideal in $\K[\bar{x},\bar{y},d,W]$ generated by the polynomials
defining $\bb_G$, by the Closure Theorem (see \cite{Cox} p. 122),
one has that $\conch_G=V(I_G\cap \K[\bar{x},d])$. Moreover,
reasoning as in Theorem 6 in \cite{SSS06} which is a direct
consequence of Exercise 7, p. 283 in \cite{Cox}, one gets that for
almost all values of $d\in \K^*$ the generic conchoid specializes
properly (see Example \ref{Example-1}). An example where the
specialization improperly behaving is taking $\cc$ as a circle
centered at $A$ and $d$ its radius (compare to Theorem
\ref{th-comp}). A similar reasoning might be done with a generic
focus, nevertheless we do not consider this situation here.
 \qed
\end{itemize}
\end{remark}

Let us illustrate the definition by two examples.

\begin{example}\label{Example-1} Let ${\cal C}$ be the  parabola over ${\Bbb C}$ defined by
$f(y_{1},y_{2})=y_{2}-y_{1}^2$, let $A=(0,-1)$ and $d=1/2$. Then
$\frak{B}(\cc,A,d)$ is defined by polynomials
\[ f(\bar{y}),\,\, C(\bar{x},\bar{y}):=(x_1-y_1)^2+(x_2-y_2)^2-\frac{1}{4}, \] \[
L(\bar{x},\bar{y}):=(y_2+1)(x_1-y_1)-y_1(x_2-y_2),\,\,
T(\bar{y},W):=W(y_{1}^{2}+(y_{2}+1)^2)-1.
\]
Now, considering $W>y_1>y_2>x_1>x_2$, and computing a Gr\"obner
basis w.r.t. the lex order, one gets that  $\conch$ is defined by
the polynomial (see Fig. \ref{fig-parabola}):

\noindent
$g(x_{1},x_{2})=16x_{1}^8+32x_{2}^2x_{1}^6+16x_{1}^4x_{2}^4+32x_{2}x_{1}^6-32x_{1}^2x_{2}^5
+24x_{1}^6-24x_{1}^4x_{2}^2-
96x_{1}^2x_{2}^4+16x_{2}^6-8x_{2}x_{1}^4-120x_{1}^2x_{2}^3+64x_{2}^5+
25x_{1}^4-68x_{1}^2x_{2}^2+92x_{2}^4+48x_{2}^3+12x_{1}^2-8x_{2}^2-16x_{2}-4,$

\noindent which is  an irreducible curve over $\Bbb C$. Similarly,
one gets that $\conch_G$ is given by

\begin{figure}[ht]
\begin{center}
\centerline{\psfig{figure=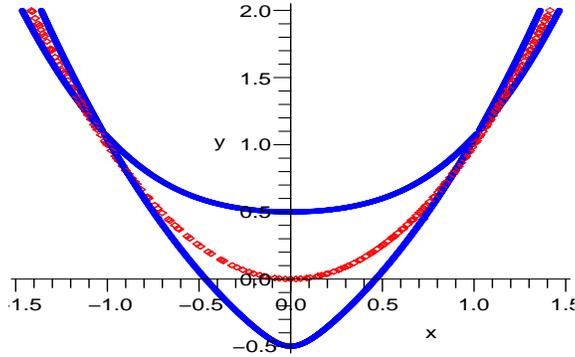,width=7.5cm,height=5.5cm}}
\end{center}
\caption{$y_{2}=y_{1}^{2}$ (in dots) and
$\mathfrak{C}(y_2-y_{1}^{2},(0,-1),1/2)$
 (in continuous traced).}
 \label{fig-parabola}
\end{figure}

\noindent $g_G(x_1,x_2,d)=-{d}^{2}+x_{2}^{2}+4\,x_{2}^{3}-2\,
x_{2}\,x_{1}^{2}- 4\,
x_{2}\,{d}^{2}-6\,x_{2}^{2}x_{1}^{2}-6\,x_{2}^{2}
{d}^{2}+6\,x_{2}^{4}+8\,x_{1}^{2}
x_{2}\,{d}^{2}-8\,x_{1}^{2}x_{2}^{3}+3\,x_{1}^{2}{d}^{2}+x_{1}^{4}+4\,
x_{2}^{5}-x_{2}^{2}x_{1}^{4}-6\,x_{2}^{4}x_{1}^
{2}+7\,x_{2}^{2}x_{1}^{2}{d}^{2}-2\,x_{1}^{4}x_{2}\,
{d}^{2}+2\,x_{1}^{6}+2\,x_{1}^{6}x_{2}+2\,x_{1}^{4}{
d}^{2}-2\,x_{1}^{2}x_{2}^{5}+2\,x_{1}^{6}x_{2}^{2}
+x_{1}^{4}x_{2}^{4}-2\,x_{2}^{2}{d}^{2}x_{1}^{4}+2
\,x_{2}^{3}{d}^{2}x_{1}^{2}+x_{2}^{6}-x_{2}^{4}{d}
^{2}-2\,x_{1}^{6}{d}^{2}+x_{1}^{8}+x_{1}^{4}{d}^{4}-4\,
x_{2}^{3}{d}^{2}$.

\noindent Note that $g_G(\bar{x},1/2)=g(\bar{x})$. \qed
\end{example}

In the following example we consider that  $\cc$ is a line which
does not pass through the focus, appearing the  well known  {\it
Conchoid of Nicomedes}.

\para

\begin{example}\label{Example-3} {\sf (Conchoid of Nicomedes)}
Let $\cc$ be the lined defined by  $f(y_{1},y_{2})=y_{2}$ and let
$A=(0,1)$. Then, $\frak{C}(\cc,A,2)$ is defined by
 (see Fig. \ref{fig-nicomedes}):
 \[
 g(x_{1},x_{2})=x_{2}^{2}x_{1}^{2}+x_{2}^{4}-2x_{2}^{3}-3x_{2}^{2}+8x_{2}-4.
 \]
\begin{figure}[ht]
\begin{center}
\centerline{\psfig{figure=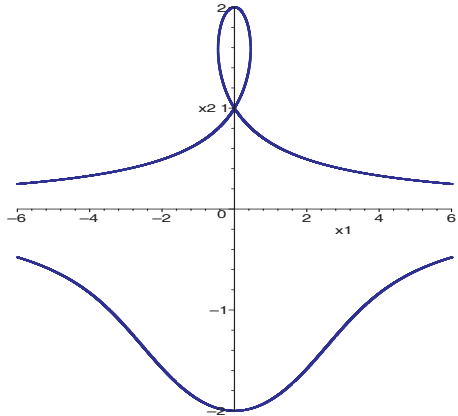,width=6cm,height=5.5cm}}
\end{center}
\caption{Conchoid of the line $y_{2}=0$ with $A=(0,1)$ and $d=2$.}
 \label{fig-nicomedes}
\end{figure}
\end{example}

We finish this section studying the connection of conchoids to
offsets (see \cite{ASS96} for further details on offsets). We
represent the offset to $\cc$ at distance $d$ as ${\cal O}_d(\cc)$.

\begin{lemma}\label{lema-conch-offset}
Let $\cc$ be a real irreducible curve, and let $A\in {\Bbb R}^2,
d\in {\Bbb R}^*$. Then, $\conch={\cal O}_{d}(\cc)$ if and only if
$\cc$ is a circle centered at $A$.
\end{lemma}
\noindent {\bf Proof.} If $\cc$ is a circle centered at $A$ the
result is obvious. Let $\conch={\cal O}_{d}(\cc)$, the result
follows from Lemma 3 in \cite{SSS05}. \qed

\section{Basic Properties of Conchoids.}\label{sec-basic-prop-conchoids}

In this section we state the first basic properties on conchoids.
For this purpose, we assume that $\cc$ is irreducible and given by
$f(y_1,y_2)$, and that $A=(a,b)$ and $d\in \K^*$ are fixed. In
addition, we consider the following Zariski open subset of  $\cc$
\[ \cc_{0}=\{ (p_1,p_2)\in \cc \,|\, (p_1-a)^2+(p_2-b)^2\neq 0\}; \]
i.e. $\cc_{0}$ consists in those points $P\in \cc$ such that $A-P$
is non-isotropic. $\cc_0$ plays an important role in the conchoid
construction. Note that if $\cc_0=\emptyset$ then the last equation
of $\bb$ does not hold, and hence $\conch$ is empty. This is the
case of the lines $(y_{1}-a)\pm \sqrt{-1}(y_{2}-b)=0$. The next
proposition relates $\cc_0$ with these lines, that we denote in the
sequel as $\lmas$ and $\lmenos$.

\begin{proposition}\label{pro-C0}
It holds that:
\begin{itemize} \item[(1)] $\cc_{0}=
\emptyset$ if and only if $\cc$ is either  $\lmas$ or $\lmenos$.
\item[(2)] $\cc \setminus \cc_{0}=\cc \cap (\lmas\cup \lmenos
)$. Moreover, $\card(\cc \setminus \cc_{0})\leq 2\,\deg(\cc)$
\item[(3)] $A\notin \cc_{0}$.
\item[(4)] If $\cc$ is real then ${\cal C}_{0}\neq \emptyset $.
\end{itemize}
\end{proposition}
\noindent   {\bf Proof.} (1) The right-left implication is trivial.
Conversely, let $\cc_0=\emptyset$. Let $g(y_{1},y_{2})=
(y_{1}-a)^2+(y_{2}-b)^2$.  Then, for every  $P\in {\cal C}$,
$g(P)=0$. That is, $g\in {\cal I}(\cc) $. Now, the proof ends taking
into account that $f$ is irreducible. (2), (3) and (4) follows from
(1). \qed

\para

In the following, we analyze the maps appearing in the incidence
diagram (see Section \ref{sec-def-conchoid}).

\begin{lemma}\label{lemma-pi}
 Let $\pi_{1}$ and $\pi_{2}$  be the projections in the incidence
 diagram associated with $\cc$.
\begin{itemize}
\item[(1)] If ${\cal C}$ is not  the circle centered at $A$ and radius $d$ then it holds that
\begin{itemize}
\item[(1.1)] $\pi_{1}$ is, at most $(2:1)$, over all points in $\pi_{1}(\bb)\setminus \{ A
\}$.
\item[(1.2)] If $A \in \pi_{1}(\bb)$ then $\card(\pi_{1}^{-1}(A))\leq
2\deg(\cc)$.
\end{itemize}
\item[(2)] $\pi_{2}$ is $(2:1)$ over all points in   ${\cal
C}_{0}$.
\end{itemize}
\end{lemma}
\noindent {\bf Proof.}
 We consider the polynomials
  \[
 \begin{array}{lcl}
C(\bar{x},\bar{y})&=&(x_{1}-y_{1})^2+ (x_{2}-y_{2})^2-d^2, \\
L(\bar{x},\bar{y})&= &(y_{2}-b) (x_{1}-y_{1})- (y_{1}-a)
(x_{2}-y_{2})  \\ &=& (b-x_{2})y_{1}+(x_{1}-a)y_{2}-bx_{1}+ax_{2}.
  \end{array}
 \]

\noindent Let us prove (1). For  $\bar{x}^{0} \in \pi_{1}(\bb)$, let
${\cal D}_{\bar{x}^{0}}$, ${\cal L}_{\bar{x}^{0}}$
be the algebraic set in $\K^2$ defined by $C(\bar{x}^{0},\bar{y})$,
$L(\bar{x}^{0},\bar{y})$,
respectively. Note that $L(\bar{x}^0,\bar{y})$ is identically zero
if and only if $\bar{x}^{0}=A$. Now, if $\bar{x}^{0}\neq A$ then
$\pi_{1}^{-1}(\bar{x}^{0})$ consists in those
$(\bar{x}^0,\bar{y},w^0)$ where $\bar{y}\in \cc \cap {\cal
D}_{\bar{x}^{0}} \cap {\cal L}_{\bar{x}^{0}}$ and
$w^0=1/((y_{1}-a)^2+(y_{2}-b)^2)$. Thus, since $\deg({\cal
D}_{\bar{x}^0})=2$  and $\deg({\cal L}_{\bar{x}_0})=1$, then
$\card(\pi_{1}^{-1}(\bar{x}^{0}))\leq 2$. If $\bar{x}^{0}=A$ then
 $\pi_{1}^{-1}(\bar{x}^{0})$ consists in those $(\bar{x}^0,\bar{y},w^0)$ where
$\bar{y}\in \cc \cap {\cal D}_{\bar{x}^{0}}$ and
$w^0=1/((y_{1}-a)^2+(y_{2}-b)^2)$. Thus, since by hypothesis $\cc$
and ${\cal D}_{\bar{x}^0}$ do not have common components, then
$\card(\pi_{1}^{-1}(\bar{x}^{0}))\leq 2\deg(\cc)$.

\noindent Now we prove (2).  Let  $\bar{y}^{0} \in \cc_{0}$, and let
${\cal D}^{\bar{y}^{0}}$ and ${\cal L}^{\bar{y}^{0}}$  be the
algebraic sets defined by $C(\bar{x},\bar{y}^{0})$, and
$L(\bar{x},\bar{y}^{0})$ respectively. ${\cal D}^{\bar{y}^{0}}$ is
the circle centered in $\bar{y}^{0}$ and radius $d$, and ${\cal
L}^{\bar{y}^{0}}$ is the line passing through $\bar{y}^{0}$ and $A$;
note that by Prop. \ref{pro-C0} $A\notin \cc_{0}$, and hence
$\bar{y}^{0}\neq A$. So, ${\cal L}^{\bar{y}^{0}}$ is not tangent to
${\cal D}^{\bar{y}^{0}}$. Moreover, $\pi_{2}^{-1}(\bar{y}^{0})$
consists in those $(\bar{x},\bar{y}^{0},w^0)$ such that $\bar{x}\in
{\cal D}^{\bar{y}^{0}} \cap {\cal L}^{\bar{y}^{0}}$ and
$w^0=1/((y^{0}_{1}-a)^2+(y^{0}_{2}-b)^2)$, with
$\bar{y}^0=(y^{0}_{1},y^{0}_{2})$. Therefore,
$\card(\pi_{1}^{-1}(A))= 2$.  \qed

\begin{remark}\label{Remark-al-lema-diagrama} Because of the last equation of $\bb$,
$\pi_{2}(\bb)\subset \cc_{0}$, and by Lemma \ref{lemma-pi} (2),
$\pi_{2}(\bb)=\cc_0$.  \qed
   \end{remark}

The following theorem essentially states that the conchoid is a
curve with at most two irreducible components.

\begin{theorem}\label{th-comp}
  Let $\cc_{0}\neq \emptyset$, then it follows that:
\begin{itemize}
\item[(1)] All the components of $\bb$  have  dimension $1$.
\item[(2)] If $\cc$ is not a circle centered at  $A$ and
radius $d$, all the components of $\conch$ have dimension $1$.
\item[(3)] If $\cc$ is  a circle centered at $A$ and
radius $d$, $\conch$ decomposes as the union of $\{A\}$ and the
circle centered at  $A$ and radius $2d$.
\item[(4)] $\conch$ has at most two components.
\end{itemize}
\end{theorem}

\noindent  {\bf Proof.} Let us prove (1). By Remark
\ref{Remark-al-lema-diagrama}, since $\cc_0\neq \emptyset$, one gets
that $\bb\neq \emptyset$.
 Let $\Gamma$ be an irreducible component of $\bb$. Let $M \in \Gamma
 $ and $P=\pi_{2}(M)\in \cc_{0}$. Let ${\cal P}(t)=(P_{1}(t), P_{2}(t))$ be a place of
 $\cc$ centered at $P$. We consider
\[ {\cal Q}^{\pm}(t)=\left({\cal P}(t)\pm
\frac{d}{\sqrt{\Delta(t)}}\cdot({\cal P}(t)-A), {\cal P}(t),
\frac{1}{\Delta(t)}\right),
\]
where $\Delta(t)=(P_{1}(t)-a)^{2}+(P_{2}(t)-b)^{2}$.
 Note that $\Delta(t)=(p_{1}-a)^{2}+(p_{2}-b)^{2}+
\cdots $, where $P=(p_{1},p_{2})$. Thus, since  $P\in \cc_{0}$, the
above power series is a unit. Therefore, each component of
 ${\cal Q}^{\pm}(t)$ can be written as a power series. Thus, ${\cal
 Q}^{\pm}(t)$ parametrize locally, respectively, two curves
 contained in $\bb$ and passing through each  of the two
 points in $\pi_{2}^{-1}(P)$. Let ${\cal Q}^{+}(t)$ be the one
 centered at $M\in  \pi_{2}^{-1}(P)$. Thus, $\dim{\Gamma}\geq 1$.
Now, let us assume that $\dim(\Gamma)> 1$, and let $\tilde{\cc}$ be
an irreducible component of $\pi_{2}(\Gamma)$. By Theorem 7 pp. 76
in \cite{Shafarevich77}, there exists an open set $U \subset
\tilde{\cc}$ such that for every $\bar{y}^{0} \in U $ it holds that
$\dim(\pi_{2}^{-1}(\bar{y}^{0}))\geq 1$. Then, taking $\bar{y}^{0}
\in U \subset \cc_{0}$ one gets a contradiction, since by Lemma
\ref{lemma-pi} (2), $\dim(\pi_{2}^{-1}(\bar{y}^{0}))=0$.

\noindent  (2) follows as (1), using that $\pi_{1}$ is always finite
(see Lemma \ref{lemma-pi} (1)).

\noindent (3) It is trivial.

\noindent For (4), the reasoning is analogous to Theorem 1 in
\cite{SS99}, using statements (1), (2) and (3) in this theorem and
Lemma \ref{lemma-pi}. \qed

\para

Next lemma follows from Lemma \ref{lemma-pi}, Theorem \ref{th-comp},
and the theorem on the dimension of fibres in \cite{Shafarevich77}
(see Theorem 7, pp.76).

\begin{lemma}\label{lema-densos-via-pi}
Let $\cc_0\neq \emptyset$ and $\pi_{1}$ and $\pi_{2}$  the
projections in the incidence diagram associated with $\cc$.
\begin{itemize}
\item[(1)] If $\Omega$ is a non-empty open subset of $\cc$, then
$\pi_{1}(\pi_{2}^{-1}(\Omega))$ is a non-empty Zariski dense subset
of $\conch$.
\item[(2)] If $\cc$ is not   the circle centered at  $A$ and radius $d$, and
$\Omega$ is a non-empty open subset of an irreducible component of
$\conch$, then $\pi_2(\pi_{1}^{-1}(\Omega))$ is a non-empty Zariski
dense subset of $\cc$.
\end{itemize}
\end{lemma}
\noindent {\bf Proof.} (1) follows using that
$\pi_{2}^{-1}(\cc_0)=\bb$ (see Remark
\ref{Remark-al-lema-diagrama}), and that $\pi_{1}(\bb)$ is
constructible in $\conch$ (see Theorem 3.16. in \cite{harris}).
\newline
\noindent In order to prove (2), let $\cal M$ be an irreducible
component of $\conch$,  let $\emptyset \neq \Omega\subset \cal M$ be
open in $\cal M$, and let $\Omega'=\Omega\cap \pi_{1}(\bb))$. Since
$\dim({\cal M})=1$ (see Theorem \ref{th-comp}(2)) and $\pi_1$ is
finite over $\Omega'$ (see Lemma \ref{lemma-pi}), at least one
component of $\pi_{1}^{-1}(\Omega')$ has dimension 1. Let $\Gamma
\subset \pi_{1}^{-1}(\Omega')$ be irreducible of dimension 1. Then,
since $\pi_2$ is finite over $\cc_0$ (see Lemma \ref{lemma-pi}), by
Theorem 7 (ii), pp. 76, in \cite{Shafarevich77}, one has that
$\dim(\pi_{2}(\Gamma))=1$. Therefore, the result follows taking into
account that $\pi_{2}(\Gamma)\subset
\pi_{2}(\pi_{1}^{-1}(\Omega'))\subset
\pi_{2}(\pi_{1}^{-1}(\Omega))\subset \cc$ and that $\cc$ is
irreducible. \qed

\begin{example}\label{Example-2} {\sf (Lima\c{c}ons of Pascal)}
 Let ${\cal C}$ be the circle centered at $(0,0)$  and radius $r=2$.
Then, the conchoid  of  ${\cal C}$ with $A=(-2,0)\in {\cal C} $ and
$d=1$ ({\it Lima\c{c}on of Pascal} at distance $1$, see Fig.
\ref{fig-caracol}), is defined by the polynomial:
 \[
 g(x_{1},x_{2})=
x_{1}^4+2x_{2}^2x_{1}^2-9x_{1}^2-4x_{1}-9x_{2}^2+12+x_{2}^4.\]

\begin{figure}[ht]
\begin{center}
\centerline{ \psfig{figure=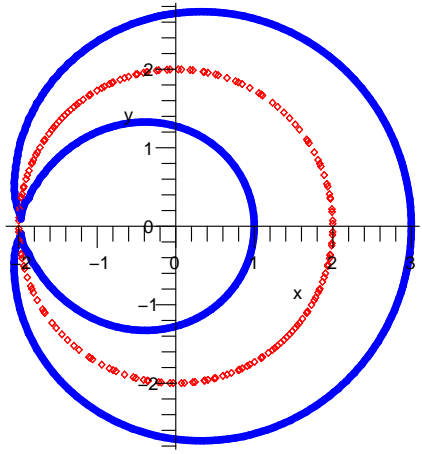,width=6cm,height=6cm}}
\end{center}
\caption{Circle centered at $(0,0)$ and radius $r=2$ (in dots), and
its conchoid from the focus $(-2,0)$ and radius $d=1$ (continuous
traced).}
 \label{fig-caracol}
\end{figure}

On the other hand, if we move the focus to  $A=(0,0)\notin {\cal
C}$, one observes that the conchoid at distance  $d=1$ has two
irreducible components (two circles centered at $A$ and radius $1$
and $3$, respectively) defined by the irreducible factors
$(x_{1}^2-9+x_{2}^2)\cdot(x_{1}^2-1+x_{2}^2).$ Note that, in this
case,  $\conch$ is the  offset of $ {\cal C}$ at distance $d=1$ (see
Lemma \ref{lema-conch-offset}). \qed
\end{example}

\section{Simple and Special Components.}\label{sec-simple-special}

In Theorem \ref{th-comp}, we have seen that if $\cc$ is not a circle
centered at the focus, all components of the conchoid are curves. In
the sequel, we assume that $\cc$ is not such a circle  centered at
the focus and that $\cc_0\neq \emptyset$. Now, we introduce and
analyze the notion of simple and special components of a conchoid.
Special and simple components provide information on the
birationality of the projections in the incidence diagram (see
Section \ref{sec-role}), and they can be used to decide whether a
curve is a conchoid (see Section \ref{sec-detecting}).  Essentially,
one component of the conchoid is special if its points are generated
for more than one point of the original curve. This phenomenon
appears when one computes conchoids of conchoids (see Theorem
\ref{th-special}). In addition, Theorem \ref{th-dist} states that,
for almost every distance and with the exception of lines passing
through the focus, all the components of the conchoid are simple.

\begin{definition}\label{def-simple}
 An irreducible component ${\cal M}$ of
$\conch$ is called {\sf simple} if there exists a non-empty Zariski
dense subset $\Omega \subset {\cal M}$ such that if $Q\in \Omega$
then $\card(\pi_{2}(\pi_{1}^{-1}(Q)))=1$. Otherwise $\cal M$ is
called {\sf special}. \qed
\end{definition}
\begin{remark}\label{remark-def-special}
An irreducible component $\cal M$ is special iff there exists a
 Zariski dense $\emptyset\neq \Omega \subset {\cal M}$ such that
for $Q\in \Omega$, $\card(\pi_{2}(\pi_{1}^{-1}(Q)))>1$. \qed
\end{remark}

\begin{proposition}\label{prop-M0}
 Let ${\cal M}$ be an  irreducible
component of
 $\conch$. Then, it holds that:
\begin{itemize}
\item[(1)] The open subset  ${\cal M}_{0} \subset {\cal M}$ is not
empty. (See Section \ref{sec-basic-prop-conchoids} for the
definition of ${\cal M}_{0}$).
\item[(2)]  ${\cal C}$ is a component of
$\mathfrak{C}({\cal M},A,d)$.
\item[(3)] Let $d'\in \K^*$ such that $d^2\neq d'^{2}$. If  ${\cal M}'$ is an  irreducible component of
  $\mathfrak{C}({\cal M},A,d')$, then ${\cal M}'$ is  a component of
$\mathfrak{C}({\cal C},A,d+d')$ or of $\mathfrak{C}({\cal
C},A,d-d')$.
\end{itemize}
\end{proposition}
\noindent {\bf Proof.} In order to prove (1), we suppose that ${\cal
M}_{0}= \emptyset$. By Proposition \ref{pro-C0} (1), ${\cal M}$ is
either $\lmas$ or $\lmenos$. Say ${\cal M}=\lmas$, similarly if
${\cal M}=\lmenos$. Then, $\Omega:=\lmas\cap \pi_{1}(\bb)\neq
\emptyset$ and dense in $\lmas$. So, there exists $P_{0}\in
\Omega\setminus \{A\}$. Let $(P_{0},Q_{0},w_{0})\in
\pi_{1}^{-1}(P_0)$; note that obviously $\pi_{1}^{-1}(P_0)\neq
\emptyset$. Now, observe that $A,P_0$ are different points on the
line $\lmas$ and, because of the third equation of $\bb$, $P_0, Q_0,
A$ are collinear. Thus, $Q_0\in \lmas$ which impossible because of
the last equation of $\bb$.
\newline
\noindent To prove (2), let $\pi_{i}$ and $\pi_{i}^{*}$ be the
projections in the incidence diagram associated with $\conch$ and
$\mathfrak{C}({\cal M},A,d)$, respectively. Let
$\cc_1:=\pi_{2}(\pi_{1}^{-1}({\cal M}_0))$. By the statement (1),
${\cal M}_0$ is a non-empty open of the irreducible component $\cal
M$, and by Lemma \ref{lema-densos-via-pi} (2), $\cc_1$ is a
non-empty Zariski dense in $\cc$; observe that we have assumed that
$\cc$ is not the circle centered at $A$ and radius $d$. Now, we see
that for every $P\in \cc_1$ there exists $Q_0\in {\cal M}_0$ and
$w_0\in \K$ such that $(P,Q_0,w_0)\in \mathfrak{B}({\cal M},A,d)$.
Indeed, since $P\in \cc_1$ there exists $Q_0\in {\cal M}_0$ and
$w_0\in \K$ such that $(Q_0,P,w_0)\in \bb$, and since $Q_0\in {\cal
M}_0$ there exists $w_1\in \K$ such that $(P,Q_0,w_1)\in
\mathfrak{B}({\cal M},A,d)$. Therefore, $\cc_1\subset
\pi_{1}^{*}(\mathfrak{B}({\cal M},A,d))$, and taking closures
$\cc\subset \mathfrak{C}({\cal M},A,d)$.
\newline
\noindent  Finally, we prove (3). Let $\pi_{i}$, $\pi_{i}^{*}$ and
$\pi_{i}^{\pm}$ be the projections in the incidence diagram
associated with $\conch$, $\mathfrak{C}({\cal M},A,d')$, and
$\mathfrak{C}(\cc,A,d\pm d')$, respectively; note that $d^2\neq
d'^{2}$ and hence $\mathfrak{C}(\cc,A,d\pm d')$ is well defined. We
consider ${\cal
M}'_{1}=\pi_{1}^{*}(\pi_{2}^{*-1}(\pi_{1}(\pi_{2}^{-1}(\cc_0))\cap
{\cal M}_0)\cap {\cal M}'_{0}$. Taking into account that $\cc_0\neq
\emptyset, {\cal M}_0\neq \emptyset, {\cal M}'_0\neq \emptyset$, by
Lemma \ref{lema-densos-via-pi}  one gets that ${\cal M}'_{1}$ is a
non-empty Zariski dense of ${\cal M}'$. Now, let $P\in {\cal
M}'_{1}$. Then, there exists $Q_0\in
\pi_{1}(\pi_{2}^{-1}(\cc_0))\cap {\cal M}_0$ and $w_0\in \K$ such
that $(P,Q_0,w_0)\in \mathfrak{B}({\cal M},A,d')$. Moreover, since
$Q_0\in \pi_{1}(\pi_{2}^{-1}(\cc_0))\cap {\cal M}_0$, there exists
$Q_1\in \cc_0$ and $w_1\in \K$ such that $(Q_0,Q_1,w_1)\in
\mathfrak{B}(\cc,A,d)$. Let us see that either $(P,Q_1,w_1)\in
\mathfrak{B}(\cc,A,d+d')$ or $(P,Q_1,w_1)\in
\mathfrak{B}(\cc,A,d-d')$. If we would prove this then it would hold
that ${\cal M}'_{1}\subset \pi_{1}^{+}(\mathfrak{B}(\cc,A,d+d'))\cup
\pi_{1}^{-}(\mathfrak{B}(\cc,A,d-d'))$. Hence, taking closures,
${\cal M}'\subset
\mathfrak{C}(\cc,A,d+d')\cup\mathfrak{C}(\cc,A,d-d')$. Therefore,
let us prove the claim; i.e. (i) $Q_1\in \cc_0$ (this is equivalent
to the first and last equation of $\mathfrak{B}(\cc,A,d\pm d')$),
(ii) $P$ is on the circle centered at $Q_1$ are radius either $d+d'$
or $d-d'$ (second equation of $\mathfrak{B}(\cc,A,d\pm d')$), (iii)
$P,Q_1,A$ are collinear (third equation of $\mathfrak{B}(\cc,A,d\pm
d')$). Indeed, we already now that $Q_1\in \cc_0$. Moreover, since
$(P,Q_0,w_0)\in \mathfrak{B}({\cal M},A,d')$ and $(Q_0,Q_1,w_1)\in
\mathfrak{B}(\cc,A,d)$, then $P,A,Q_0$ are collinear and $Q_0,Q_1,A$
are also collinear. Therefore, since $A\not\in\{P,Q_0,Q_1\}$ (this
follows from Prop. \ref{pro-C0} (2), because $P\in {\cal
M}'_{0},Q_1\in \cc_0,$ and $Q_0\in {\cal M}_0$), then $P,Q_0,Q_1,A$
are collinear; in particular (iii) holds. Also, since
$(P,Q_0,w_0)\in \mathfrak{B}({\cal M},A,d')$ and $(Q_0,Q_1,w_1)\in
\mathfrak{B}(\cc,A,d)$, then $P$ is on the circle centered at $Q_0$
radius $d'$ and $Q_1$ is on the circle centered at $Q_0$ and radius
$d$. Thus, since $P,Q_0,Q_1$ are collinear (see above), then (ii)
holds. \qed

\begin{remark}\label{remark-cfa}
Note that, by Proposition \ref{prop-M0} (2) and Lemma
\ref{lema-conch-offset}, if $\cc$ is not a circle centered at the
focus, then for every $d\in \K^*$ none component of $\conch$ is a
circle centered at $A$. \qed
\end{remark}

\para

Next theorem shows that, similarly as in the offsetting construction
(see \cite{SS99}), special components appear only when computing
conchoids of conchoids.

\para

\begin{theorem}\label{th-special}  An  irreducible component ${\cal
M}$ of
 $\conch$ is special if an only if $\mathfrak{C}({\cal
 M},A,d)={\cal C}$.
\end{theorem}
\noindent  {\bf Proof.} Let $\cal M$ be special. We assume w.l.o.g.
that ${\cal M}_0$ is the Zariski dense where the cardinality of the
fiber is bigger than 1. By Prop. \ref{prop-M0} (2),
  ${\cal C} \subset \mathfrak{C}({\cal
 M},A,d)$. In order to see that  $\mathfrak{C}({\cal
 M},A,d)\subset \cc$, let $\pi_i$ and $\pi_{i}^*$ the projections in
 the incidence diagram associated with $\cc$ and $\cal M$,
 respectively. Let
 $\Omega^*=\pi_{1}^{*}(\pi_{2}^{*-1}({\cal M}_0))$. By Lemma \ref{lema-densos-via-pi},
  $\Omega^*$ is dense in $\mathfrak{C}({\cal M},A,d)$. Now, let $P\in \Omega^*$. Then, there exists $Q\in
 {\cal M}_0$ and $w_0\in \K$ such that $(P,Q,w_0)\in \mathfrak{B}({\cal
 M},A,d)$. Moreover, using that $Q\in {\cal M}_0$,  and that $\cal M$ is special, one gets that
 $\card(\pi_{2}(\pi_{1}^{-1}(Q)))>1$. Let $P_1,P_2\in
 \pi_{2}(\pi_{1}^{-1}(Q))\subset \cc_0$. Then, there exists
 $w_1,w_2\in \K$ such that $(Q,P_1,w_1),(Q,P_2,w_2)\in \bb$.
 Therefore, $Q,P_1,A$ are collinear, $Q,P_2,A$ are collinear, and
 $P_1,P_2$ are on the circle $\cal D$ centered at $Q$ and radius $d$.
 Furthermore, from $(P,Q,w_0)\in \mathfrak{B}({\cal
 M},A,d)$, one gets that $P,Q,A$ are collinear and $P\in \cal D$.
 Since $Q\neq A$, because $Q\in {\cal M}_0$ (see Prop.
 \ref{pro-C0}), $P_1,P_2,P,Q,A$ are on the same line $\cal L$, and
 $\{P_1,P_2,P\}\subset {\cal L}\cap {\cal D}$. Since the center of
 $\cal D$ is on $\cal L$, then $P=P_1$ or $P=P_2$. So, $P\in \cc$.
 Finally, since $\Omega^*$ is dense in $\mathfrak{C}({\cal
 M},A,d)$, taking closures one gets $\mathfrak{C}({\cal
 M},A,d)\subset \cc$.

Conversely, let $\mathfrak{C}({\cal
 M},A,d)={\cal C}$. Let $\pi_i$ and $\pi_{i}^{*}$ as above. Let
 ${\cal M}_1=\pi_{2}^{*}(\pi_{1}^{*-1}(\cc_0))\cap {\cal M}_0$. Since $\mathfrak{C}({\cal
 M},A,d) =\cc$, by Lemma \ref{lema-densos-via-pi}, ${\cal M}_1$ is a non-empty Zariski dense  in $\cal M$.  We prove
 that for every $Q\in {\cal M}_1$,
 $\card(\pi_{2}(\pi_{1}^{-1}(Q)))>1$ and hence that $\cal M$ is
 special. By Lemma \ref{lemma-pi}(2),
 $\card(\pi_{1}^{*}(\pi_{2}^{*-1}(Q)))=2$. Let $P_1,P_2\in
 \pi_{1}^{*}(\pi_{2}^{*-1}(Q))\subset  \cc_0$, $P_1\neq P_2$. Then, there exists $w_0\in \K$ such that
  $(P_1,Q,w_0),(P_2,Q,w_0)\in \mathfrak{B}({\cal
 M},A,d)$. So, $P_1,P_2,Q$ satisfy the second and third
 equation of $\bb$. Moreover, since $P_i\in\cc_0$, there exists $w_i\in \K$ such that $(Q,P_i,w_i)\in
 \bb$, and therefore $P_1,P_2\in \pi_{2}(\pi_{1}^{-1}(Q))$. \qed

\para

We illustrate the previous results by an example.

\begin{example}\label{Example-4}
Let $\cc$ be the circle centered at $(0,0)$ and radius $r=1$ defined
by the polynomial $f(y_{1},y_{2})=y_{1}^2+y_{2}^2-1$. Let
$A=(-1,0)\in \cc$. First, we compute  $\frak{C}(\cc,(-1,0),2)$,
obtaining  a Lima\c{c}on of Pascal, defined by the polynomial:
\[
g(x_{1},x_{2})=x_{1}^4+2\,x_{2}^2\,x_{1}^2-6\,x_{1}^2-8\,x_{1}-6\,x_{2}^2-3+x_{2}^4.
\]
Note that we get  a cardioid. Let us denote it as $\cc^{'}$, i.e.
$\cc^{'}=\frak{C}(\cc,A,2)$ (see Fig. \ref{fig-cardioide} left).

\begin{figure}[ht]
\begin{center}
\centerline{ \psfig{figure=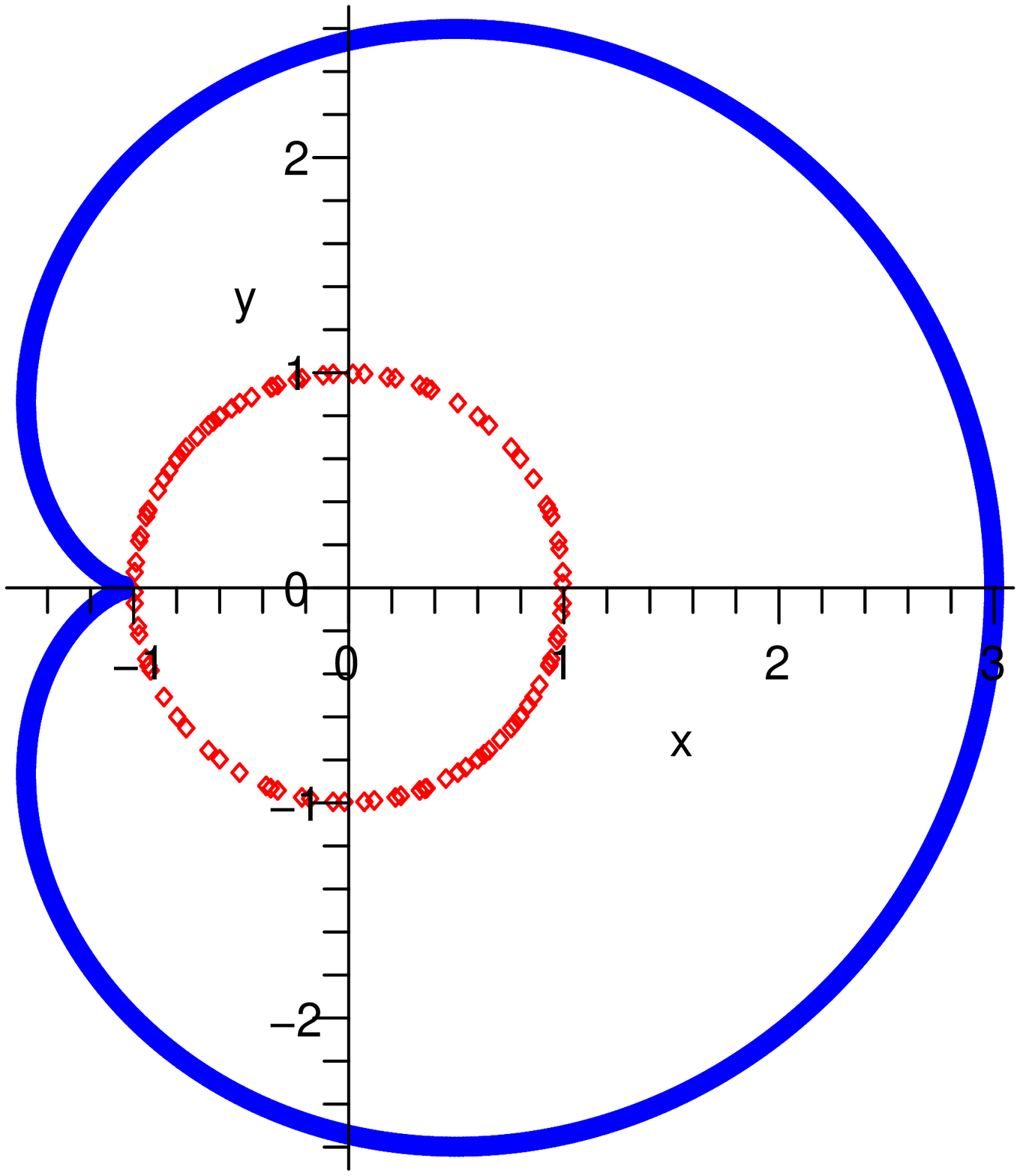,width=4.5cm,height=4.5cm}
\psfig{figure=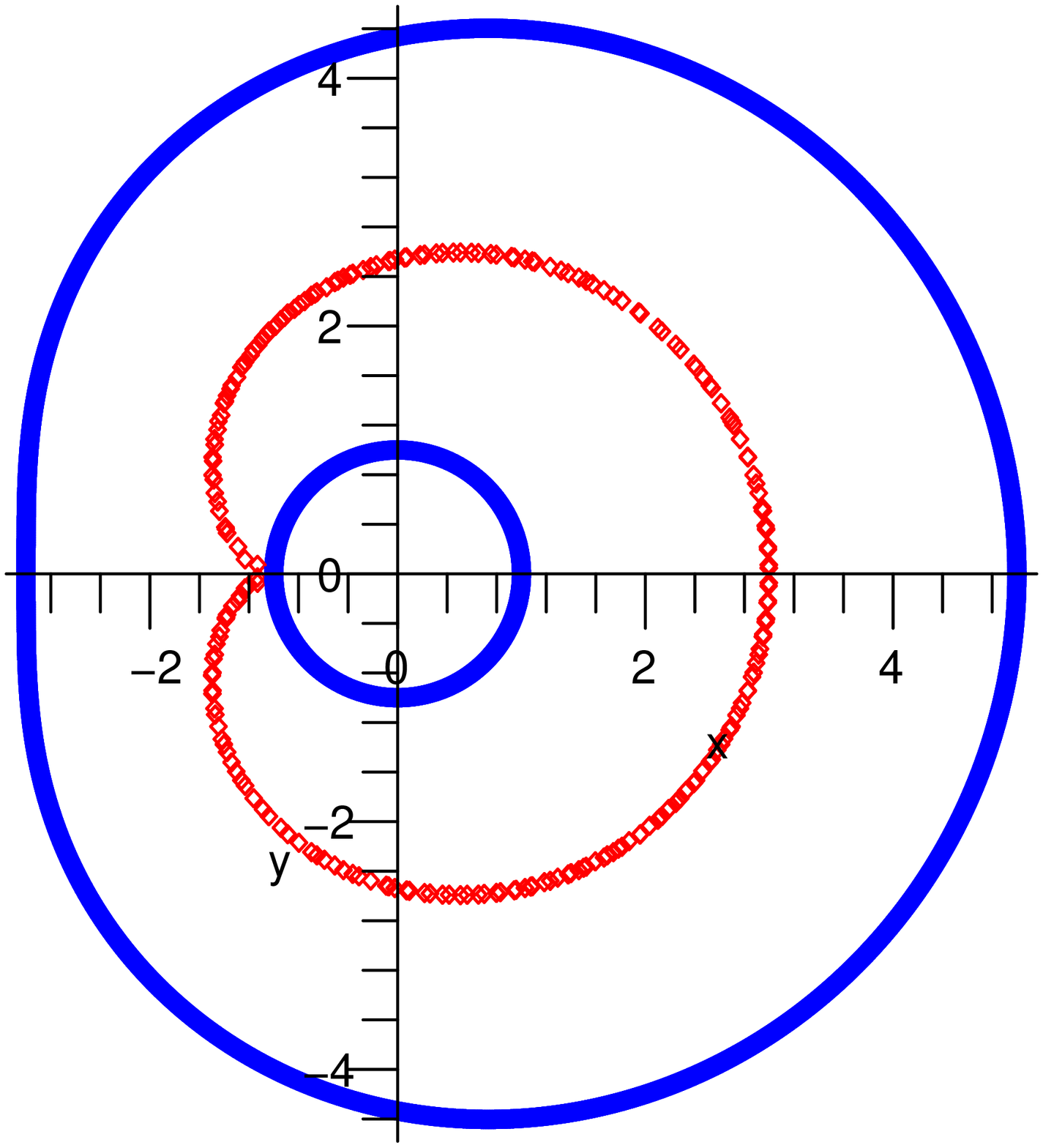,width=4.5cm,height=4.5cm}
\psfig{figure=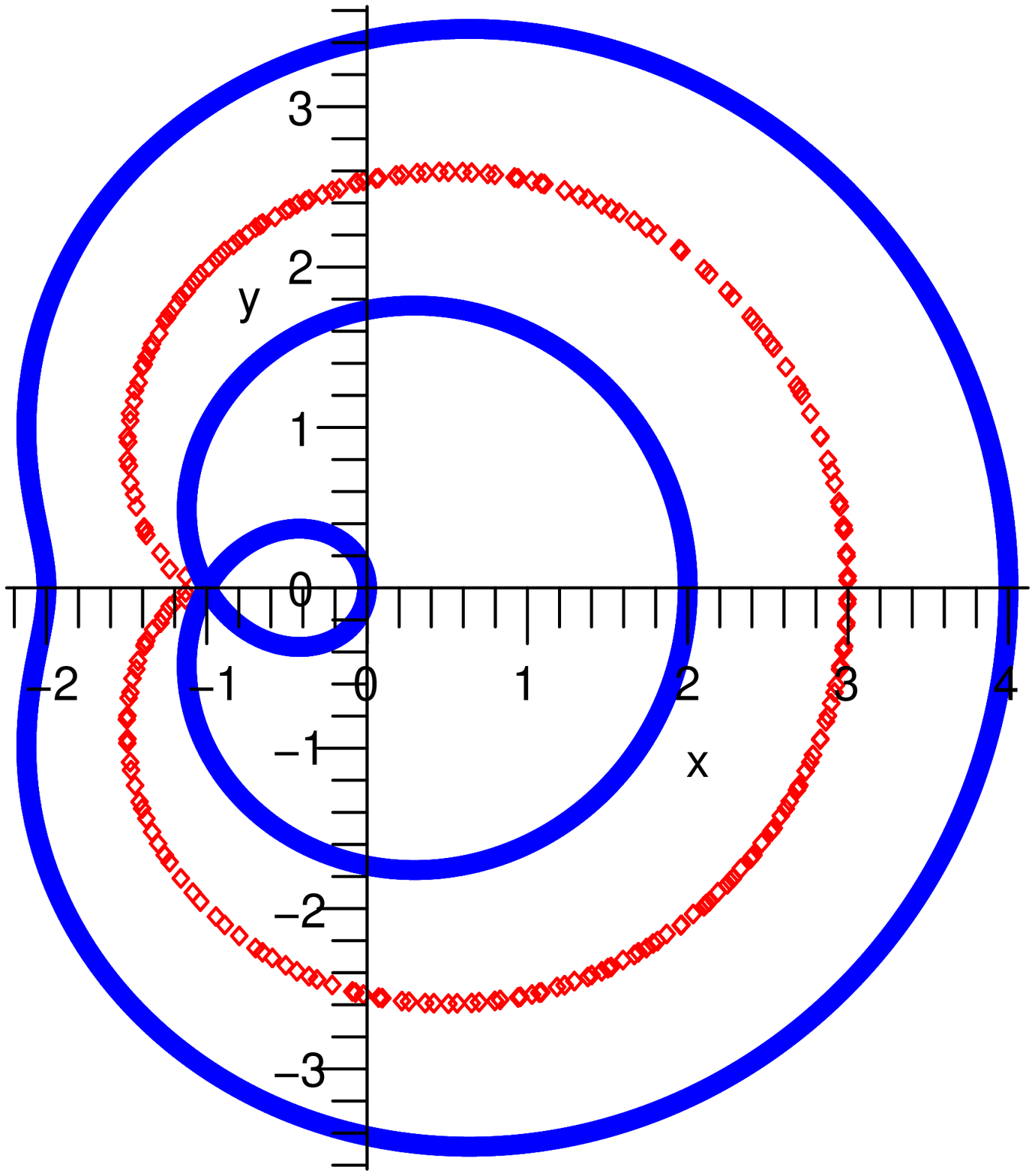,width=4.5cm,height=4.5cm}}
\end{center}
\caption{{\sf Left:} Circle centered at $(0,0)$ and radius $r=1$ (in
dots), and its Conchoid from the focus $(-1,0)$ and radius $d=2$
(cardioid). \newline {\sf Center:} Cardioid (in dots), and its
Conchoid from the focus $(-1,0)$ and radius $d=2$ (continuous
traced). \newline {\sf Right:} Cardioid (in dots), and its Conchoid
from the focus $(-1,0)$ and radius $d=1$ (continuous traced). }
 \label{fig-cardioide}
\end{figure}

Now, we compute the conchoid of the cardioid $\cc{'}$ from the same
focus $A=(-1,0)$ and the same distance $d=2$.   In this case, one
gets a reducible curve, say $\cc^{''}:=\frak{C}(\cc^{'},A,2)$, with
two irreducible components defined by the irreducible factors of:
\[(x_{1}^2+x_{2}^2-1)\cdot(x_{1}^4+2\,x_{2}^2\,x_{1}^2-18\,x_{1}^2-32\,x_{1}-18\,x_{2}^2-15+x_{2}^4).
\]
Note that, one component of $\cc^{''}$ is the initial circle $\cc$
(and therefore it is an special component of $\cc^{''}$; see Theorem
\ref{th-special}), and the other component is a Lima\c{c}on of
Pascal of $\cc$ from the focus $A$ and distance $d=4$. See Fig.
\ref{fig-cardioide} center.


On the other hand, the conchoid of $\cc^{'}$, from the focus $A$ but
now taking distance $d=1$, decomposes as the union of two
irreducible  components defined by the irreducible factors of the
equation:
\[(x_{1}^4+2\,x_{2}^2\,x_{1}^2-11\,x_{1}^2-18\,x_{1}-11\,x_{2}^2-8+x_{2}^4)\cdot
(x_{1}^4+2\,x_{2}^2\,x_{1}^2-3\,x_{1}^2-2\,x_{1}-3\,x_{2}^2+x_{2}^4),
\]
that  correspond to two Lima\c{c}ons of Pascal of $\cc$ from the
focus $A$ at distance $d=3$ and $d=1$ respectively, as indicates
 Proposition \ref{prop-M0} (3); see Fig. \ref{fig-cardioide} right. \qed

\end{example}
\para

\para

Next theorem states the main property of the components of a
conchoid.

\para

\begin{theorem}\label{th-simple} Let ${\cal C}$ be different to
a line passing through the focus, then  $\conch$ has at least one
simple component.
\end{theorem}
\noindent  {\bf Proof.} Recall that we have assumed that  $\cc$ is
not a circle centered at $A$.  By Theorem \ref{th-comp} (4),
$\conch$ has at most two irreducible components. So, we distinguish
two cases: (i) $\conch={\cal M}$ and (ii) $\conch={\cal M}'\cup
{\cal M}''$; where all components are taken irreducible. Let us
assume that all components of $\conch$ are special. By Theorem
\ref{th-special} one has that (i) $\mathfrak{C}({\cal M},A,d)={\cal
C}$ or (ii) $\mathfrak{C}({\cal M}',A,d)=\mathfrak{C}({\cal
M}'',A,d)={\cal C}$. Let $r=\deg(\cc)$.
 We take  $P\in
 \cc_0$ such that $\Delta:=+\|P-A\|\not\in \{-n\,d \,|\, n\in {\Bbb
 N}\,\,\mbox{and}\,\, 1\leq n\leq 2r+1\}$. Note that this is always
 possible because $\cc_0\neq \emptyset$ and $\cc$ is not a circle
 centered at $A$. Also, let $\pi_i$  be the projections in the
incidence diagram of
 $\cc$, and let $\pi_{i}^{*}$ be the projections in the incidence
 diagram of $\mathfrak{C}({\cal M},A,d)$ (if (i)
 happens) and $\pi_{i}',\pi''_{i}$ be the projections in the incidence diagrams of  $\mathfrak{C}({\cal M}',A,d)$ and
 $\mathfrak{C}({\cal M}'',A,d)$, respectively (if (ii)
 happens). We consider the families of points
\[
\{P_n:=P+\frac{2n\,d}{\Delta}(P-A)\}_{n\in {\Bbb N}, 0\leq n \leq
r},\,\, \{Q_n:=P+\frac{(2n+1)\,d}{\Delta}(P-A)\}_{n\in {\Bbb N},
0\leq n \leq r}. \] We prove by induction that $P_n\in \cc_0$ and
$Q_n\in \conch_0$. Observe that, if this holds, then the line
passing through $A$ and $P$ intersects $\cc$ at $r+1$ different
points, which is impossible because $\cc$ is irreducible and it is
not a line passing though $A$. Let us prove the induction claim.
$P_0=P\in \cc_0$. Furthermore, by construction $P_0,Q_0, A$ are
collinear and $\|Q_0-P_0\|=d$. Therefore, $(Q_0,P_0,1/\Delta)\in
\bb$, and hence $Q_0=\pi_{1}(Q_0,P_0,1/\Delta)\in \conch$. Moreover,
$\|Q_0-A\|^2=(\Delta+d)^2$, that is not zero because of the
selection of $P$ (see above), and hence $Q_0\in \conch_0$. Now, we
prove it for $P_n$ and $Q_n$. Since $P_n=Q_{n-1}+d/\Delta\, (P-A)$,
then $P_n, Q_{n-1}, A$ are collinear, and $\|P_{n}-Q_{n-1}\|=d$.
Moreover, by induction hypothesis, $Q_{n-1}\in\conch_0$. Now, if (i)
happens, then $Q_{n-1}\in {\cal M}_0$. Thus, there exists $w\in \K$
such that $(P_n,Q_{n-1},w)\in \mathfrak{B}({\cal M},A,d)$. So,
$P_n=\pi_{1}^{*}(P_n,Q_{n-1},w)\in \mathfrak{C}({\cal M},A,d)=\cc$.
In addition, $\|P-A\|^2=(\Delta+2nd)^2\neq 0$, because of the
selection of $P$. Hence $P_n\in \cc_0$. On the other hand, if (ii)
happens, then $Q_{n-1}\in {\cal M}'_0\cup {\cal M}''_{0}$. Say
$Q_{n-1}\in {\cal M}'_{0}$, similarly if $Q_{n-1}\in {\cal
M}''_{0}$. Thus, there exists $w\in \K$ such that
$(P_n,Q_{n-1},w)\in \mathfrak{B}({\cal M}',A,d)$. So,
$P_n=\pi'_{1}(P_n,Q_{n-1},w)\in \mathfrak{C}({\cal M}',A,d)=\cc$. In
addition, reasoning as above $\|P-A\|^2=(\Delta+2nd)^2\neq 0$. Hence
$P_n\in \cc_0$. Furthermore, since $Q_n=P_n+d/\Delta\,(P-A)$, then
$P_n,Q_n,A$ are collinear and $\|Q_n-P_n\|=d$. Since, we already
know that $P_n\in \cc_0$, there exists $w\in \K$ such that
$(Q_n,P_n,w)\in \bb$. So, $Q_n=\pi_{1}(Q_n,P_n,w)\in \conch$.
Finally, $\|Q_n-A\|^2=(\Delta+(2n+1)d)^2\neq 0$, because of the
selection of $P$. Thus, $Q_n\in \conch_0$. \qed

\para

The next corollary follows from the previous theorem and Theorem
\ref{th-comp}.

\begin{corollary}\label{corolario-simple} Let ${\cal C}$ be different to
a line passing through the focus, then it holds that:
\begin{itemize}
\item[(1)] If $\conch$ is irreducible
then is  simple.
\item[(2)] $\conch$ has at least one  simple component.
\end{itemize}
\end{corollary}

\begin{remark}\label{Remark-3} Let $\cc$ be a line different to $\lmas$ and $\lmenos$ (i.e. $\cc_0\neq \emptyset$), and let $A\in \cc$.
Then $\conch$ is irreducible (in fact $\conch=\cc$) and special for
every $d\neq 0$. Indeed, let $\cc$ be the line of equation
\[ f=\lambda (y_1-a)+\mu (y_2-b).  \]
Observe that, by Prop. \ref{pro-C0}, $\cc_0=\cc\setminus \{A\}$.
Then, if $\lambda=0$, one gets that
\[ \pi_{1}(\bb)=\{(y_1\pm d, b)\,|\, (y_1,y_2)\in
\cc\setminus \{A\}\} \subset \cc, \] and hence $\conch=\cc$.
Moreover, if $Q=(x_1,x_2)\in \pi_{1}(\bb)\setminus \{(a\pm d,b)\}$
then $\pi_{2}(\pi_{1}^{-1}(Q))=\{(x_1-d,b),(x_1+d,b)\}$. Thus
$\conch$ is special. On the other hand, if
 $\lambda\neq 0$, one gets that
\[ \pi_{1}(\bb)=\left\{\left(y_1\mp
\frac{d\mu}{\sqrt{\lambda^2+\mu^2}}, y_2\pm
\frac{d\lambda}{\sqrt{\lambda^2+\mu^2}}\right)\,|\, (y_1,y_2)\in
\cc_0 \right\}\subset \cc. \] Note that $\lambda^2+\mu^2\neq 0$
because $\cc\neq \lmas $ and $\cc\neq \lmenos$. Therefore,
$\conch=\cc$. Moreover, if
\[Q=(x_1,x_2)\in \pi_{1}(\bb)\setminus \left\{\left(a\mp
\frac{d\mu}{\sqrt{\lambda^2+\mu^2}}, b\pm
\frac{d\lambda}{\sqrt{\lambda^2+\mu^2}}\right)\right\}\] then
\[ \pi_{2}(\pi_{1}^{-1}(Q))=\left\{\left(x_1\pm
\frac{d\mu}{\sqrt{\lambda^2+\mu^2}}, x_2\mp
\frac{d\lambda}{\sqrt{\lambda^2+\mu^2}}\right)\right\}\]
 and hence $\conch$ is
special. \qed
\end{remark}

\begin{corollary}\label{co-lines}\para
 The only curves for which the conchoid
is irreducible and special are the lines passing through the focus.
\end{corollary}
\noindent   {\bf Proof.} It follows from Theorem \ref{th-simple} and
the above Remark. \qed

\begin{lemma}\label{lema-special} Let ${\cal C}$ be different to
a line passing through the focus, then there exist, at most, a
finite number of distances for which all  the conchoids $\conch$
have one common component.
\end{lemma}
\noindent  {\bf Proof.} Let  $D\subset \K^*$ be an infinite set  and
 $\cal M$ an irreducible algebraic curve  such that ${\cal M}\subset \bigcap_{d\in D}
\conch$. Let  $r=\deg(\cc)$, and  $d_{1},\ldots, d_{r+1} \in D$ such
that $d_{i}^{2}\neq d_{j}^{2}, \,\, \forall i\neq j$. Let $
\Omega=\bigcap_{i=1}^{r+1} \pi_{1,i}(\pi_{2,i}^{-1}(\cc_{0}))\cap
{\cal M}_{0}$, where $\pi_{1,i}$ and $\pi_{2,i}$ are the projections
of the incidence diagram associated to $\frak{C}(\cc,A,d_i)$.  Since
$\cal M$ is irreducible, since is ${\cal M}_0$ non-empty an open in
$\cal M$,  by Lemma \ref{lema-densos-via-pi} one has that $\Omega $
is a non-empty Zariski subset of $\cal M$. Now, let $Q\in \Omega$.
Note that since $Q\in {\cal M}_0$ then  $Q\neq A$. Then,  for all $i
\in \{1,\ldots, r+1\}$ there exist $P_{i}\in \cc_{0} $ and $w_{i}\in
\K^*$ such that $(Q,P_{i},w_{i})\in
 \mathfrak{B}(\cc,A,d_{i})$.
 So,   $Q, P_{i},A$ are  collinear and
  $\| Q- P_{i}\|^2=d_{i}^{2}$. Furthermore, t  $P_{i}\neq P_{j}$ for all $i\neq j$ and $i,j \in
\{1,\ldots r+1\}$, since $\| Q- P_{i}\|^2=d_{i}^{2} \neq
d_{j}^{2}=\| Q- P_{j}\|^2$. Now, let ${\cal L}$ denote the
  line passing through
  $Q$ and $A$; note that ${\cal L}$ is well-defined since $Q\neq A$.
  Then, one has that $
\{P_{1},\ldots,P_{r+1}\}\subset {\cal L}\cap \cc$. Thus,
$\card({\cal L}\cap \cc)>\deg(\cc)$ and hence $\cc={\cal L}$, in
contradiction with the hypothesis. \qed

\para

\begin{theorem}\label{th-dist} Let ${\cal C}$ be different to
a line passing through the focus. Then for almost every distance $d
\in \K^* $ all the components of the conchoid $\conch$  are simple.
\end{theorem}

\noindent  {\bf Proof.} Let us assume that there exists an infinite
set $D\subset \K^*$ such that for all $d \in \K^* $, the conchoid
$\conch$ has an special component. Let $r=\deg( \cc)$,
$\delta_r:=1+(^{r}_{2})$, and $d_{1},\ldots,d_{\delta_r} \in D$ such
that: $d_{i}^{2}\neq d_{j}^{2}$,
  $\cc$ is not a component of $\mathfrak{C}(\cc,A,d_{i})$, and
${\cal M}_{i}\neq {\cal M}_{j}$ being ${\cal M}_{i}$  the special
component of $\mathfrak{C}(\cc,A,d_{i})$. Note that this is always
possible because $d_{i}^{2}\neq d_{j}^{2}$ and  lemma
\ref{lema-special}. Let  ${\cal M}_{i,0}=\{P\in {\cal M}_i\,|\,
\|P-A\|^2\neq 0\}$, and let  \[ \Delta_{i}= ({\cal M}_i \setminus
{\cal M}_{i,0})\cup ({\cal M}_i \cap \cc)\bigcup_{j\neq i}({\cal
M}_i\cap {\cal M}_{j}). \]  Since ${\cal M}_i\neq {\cal M}_j$,
${\cal M}_i\neq \cc$,  using that
  ${\cal M}_{i,0}$ is open an non-empty in
${\cal M}_i$ (see Proposition \ref{prop-M0}), and that $\cc_0$ is
open and non-empty in $\cc$,   one has that
$\Delta:=\bigcup_{i=1}^{\delta_r} \Delta_{i}\cup (\cc\setminus
\cc_0)$ is a finite set. In this situation, we take a line $\cal L$
passing through $A$ and such that ${\cal L} \cap \Delta\subset
\{A\}$. Let
 $Q_i\in {\cal L}\cap {\cal M}_{i,0}$ for $i\in\{1,\ldots,\delta_r\}$. This, in particular, implies that $Q_i\neq A$. By construction
$\card(\{Q_1,\ldots,Q_{\delta_r}\})=\delta_r$. Since ${\cal M}_i$ is
special, by Theorem \ref{th-special} one gets that $\frak{C}({\cal
M}_i,A,d_i)=\cc$. Furthermore, by construction $Q_i\in {\cal
M}_{i,0}$. Therefore, there exist $P'_{i},P_i\in \cc$, $P'_{i}\neq
P_i$, and $w_{i}$ such that $(P_i,Q_i,w_{i}),(P'_{i},Q_i,w_{i})\in
\frak{B}({\cal M}_i,A,d_i)$. Thus,
$\|Q_i-P'_{i}\|^2=\|Q_{i}-P_{i}\|^2=d_{i}^{2}$ and
$Q_i,P_i,P'_{i}\in {\cal L}$. Let $C_i$ be the circle centered at
$Q_i$ and radius $d_i$. Since $Q_i\neq Q_j$  we have $\delta_r$
different circles with all the centers at $\cal L$. This implies
that $\card({\cal L}\bigcap_{i=1}^{\delta_r} C_i)\geq r+1$. On the
other hand ${\cal L}\cap C_i=\{P_i,P'_{i}\}\subset \cc$. Hence there
exist, at least, $r+1$ different points in ${\cal L}\cap \cc$ and
therefore $\cc ={\cal L}$, in contradiction with the hypothesis.
\qed

\section{The Role of Simple Components}\label{sec-role}
Simple components (see previous section) play an important role,
from the theoretical point of view, in the study of conchoids.
Essentially, they provide information on the birationality of the
maps in the incidence diagram, and hence they open the door for
studying, in further research, algebraic and geometric properties of
the conchoids.

\begin{lemma}\label{lemma-pi-birracional}
 Let $\pi_{1}$ and $\pi_{2}$  be the projections in the incidence
 diagram associated with $\cc$.
\begin{itemize}
\item[(1)] Let $\cc$ be different to a circle centered at $A$ and radius $d$. If $\conch$ is reducible and
${\cal M}$ is an irreducible component of $\conch$, then the
restricted map
\[\tilde{\pi}_{2}=\pi_{2}\vert _{\pi_{1}^{-1}({\cal M})}:
\pi_{1}^{-1}({\cal M})
 \longrightarrow \cc\] is birational.
\item[(2)] If ${\cal M}$ is an irreducible component of $\conch$,
then the restricted map
\[\tilde{\pi}_{1}=\pi_{1}\vert_{\pi_{1}^{-1}({\cal M})}:
\pi_{1}^{-1}({\cal M})\longrightarrow {\cal M}\] is birational iff
${\cal M}$ is simple.
\end{itemize}
\end{lemma}
\noindent {\bf Proof.} (1) Let ${\cal M}'$ be the other component of
$\conch$ (see Theorem \ref{th-comp}).    Let $\Sigma={\cal M}\cap
{\cal M}'$. By Lemma \ref{lema-densos-via-pi},
$\cc_1:=\pi_{2}(\pi_{1}^{-1}({\cal M}\setminus \Sigma)) \cap
\pi_{2}(\pi_{1}^{-1}({\cal M}'\setminus \Sigma))\cap \cc_0$ is a
non-empty Zariski dense of $\cc$. We prove that for $P\in \cc_1$ the
cardinality of $\tilde{\pi}_{2}^{-1}(P)$ is 1. Since $P\in \cc_0$,
by Lemma \ref{lemma-pi}, $\card(\pi_{2}^{-1}(P))=2$. Moreover, since
$P\in \pi_{2}(\pi_{1}^{-1}({\cal M}\setminus \Sigma)) \cap
\pi_{2}(\pi_{1}^{-1}({\cal M}'\setminus \Sigma))$, there exist $H\in
\pi_{1}^{-1}({\cal M}\setminus \Sigma))$ and $H'\in
\pi_{1}^{-1}({\cal M}'\setminus \Sigma))$ such that $\pi_{2}(H)=P$
and $\pi_{2}(H')=P$; i.e. $\pi_{2}^{-1}(P)=\{H,H'\}$. Thus, there
exist $Q\in {\cal M}\setminus \Sigma$, $Q'\in {\cal M}'\setminus
\Sigma$ (note that this implies that $Q\in {\cal M}$ but $Q'\not\in
{\cal M}$) such that $H=(Q,P,w)$ and $H'=(Q',P,w)$. So, $H\in
\pi_{1}^{-1}({\cal M})$ and $H'\not\in \pi_{1}^{-1}({\cal M})$.
Therefore $\tilde{\pi}_{2}^{-1}(P)=\{H\}$.

\noindent (2) follows  from the notions of simple component and
birational map. \qed

\para

From this lemma, one directly deduces the following corollary.

\begin{corollary}\label{cor-comp-simples-birracionales}
Let $\cc$ be such that $\conch$ is reducible. Then, the simple
components of $\conch$ are birationally equivalent to $\cc$.
\end{corollary}

In the following example we illustrate these results.

\begin{example} Let $\cc$ be the plane curve defined by
\[ f(y_1,y_2)=-3+9\,{y_{{1}}}^{2}+9\,{y_{{2}}}^{2}+2\,y_{{2}}-4\,{y_{{2}}}^{4}-4\,{y
_{{1}}}^{4}-8\,{y_{{1}}}^{2}{y_{{2}}}^{2} .\] Let $A=(0,-1)$ and
distance $d=1/2$ then $\conch={\cal M}\cup {\cal N}$ where $\cal N$
is defined by $N(\bar{x}):=x_{1}^{2}+x_{2}^{2}-1$ and $\cal M$ by
$M(\bar{x}):={x_{{1}}}^{4}+2\,{x_{{2}}}^{2}{x_{{1}}}^{2}-3\,{x_{{1}}}^{2}-3\,{x_{{2
}}}^{2}-2\,x_{{2}}+{x_{{2}}}^{4}$. Note that $\cal N$ is special
(i.e. $\frak{C}({\cal N},A,d)=\cc$) and $\cal M$ is simple. $\bb$
decomposes as $\bb=\Gamma_1\cup \Gamma_2$ where
$\Gamma_{1}:=\overline{\pi_{1}^{-1}({\cal N})}$ and
$\Gamma_2:=\overline{\pi_{1}^{-1}({\cal M})}$. By Lemma
\ref{lemma-pi-birracional} (2), $\pi_{1}\vert_{\Gamma_2}:
\Gamma_2\longrightarrow {\cal M}$ is birational and
$\pi_{1}\vert_{\Gamma_2}: \Gamma_1\longrightarrow {\cal N}$ is not.
Indeed
\[ \begin{array}{lccc}
(\pi_{1}\vert_{\Gamma_2})^{-1}: & {\cal M} & \longrightarrow &
\Gamma_2
\\ & (x_1,x_2) & \longmapsto & \left(x_{{1}},x_{{2}},-{\frac
{x_{{1}} \left( {x_{{1}}}^{2}-6+{x_{{2}} }^{2}-4\,x_{{2}} \right)
}{4(x_{{2}}+1)}},
-\frac{{x_{{2}}}^{2}-4x_{{2}}+{x_{{1}}}^{2}-2}{4},\frac{4}{8\,x_{{2}}+9}\right)
\end{array}.\] Thus, by Lemma \ref{lemma-pi-birracional} (1),
\[ \begin{array}{lccc}
\varphi: & {\cal M} & \longrightarrow & \cc
\\ & (x_1,x_2) & \longmapsto & \displaystyle{\left(-{\frac
{x_{{1}} \left( {x_{{1}}}^{2}-6+{x_{{2}} }^{2}-4\,x_{{2}} \right)
}{4(x_{{2}}+1)}},
-\frac{{x_{{2}}}^{2}-4x_{{2}}+{x_{{1}}}^{2}-2}{4}\right)}
\end{array}\]
is birational (i.e. $\varphi= \pi_{2}\vert_{\Gamma_2} \circ
(\pi_{1}\vert_{\Gamma_2})^{-1}$). In fact,
\[ \begin{array}{lccc} \varphi^{-1}: & {\cal C} & \longrightarrow
& {\cal M}
\\ & (y_1,y_2) & \longmapsto & \displaystyle{\left({\frac {y_{{1}} \left( 3+8\,y_{{2}}+4\,{y_{{1}}}^{2}+4\,{y_{{2}}
}^{2} \right)
}{8(y_{{2}}+1)}},y_{{2}}-\frac{5}{8}+\frac{1}{2}\,{y_{{1}}}^{2}+\frac{1}{2}\,{y_{{2
}}}^{2}\right)}.
\end{array}\] \qed
\end{example}
\section{Detecting Conchoids}\label{sec-detecting}
In this section, we show how special components can be used to
decide whether a given irreducible plane curve $\cal D$ (being
different to $\lmas$ and $\lmenos$) is the conchoid of another
curve. First observe that one can always find a curve $\cc$, and
$A\in \K^2$, $d\in \K^*$ such that $\cal D$ is a component of
$\conch$. For this purpose, one simply has to take $\cc$ as an
irreducible component of $\frak{C}({\cal D},A,d)$. Then, by
Proposition \ref{prop-M0} (2), ${\cal D}\subset \conch$. So, we are
now interested in deciding whether there exist $A\in \K^2$, $d\in
\K^*$ and $\cc$ such that ${\cal D}=\conch$. By Theorem
\ref{th-special}, this is equivalent to decide whether there exist
$A\in \K^2$, $d\in \K^*$ such that $\frak{C}({\cal D},A,d)$ has an
special component; if so, $\cc$ is the special component. We proceed
as follows

\para

\noindent \underline{\sf Finding the focus.} Let $g(\bar{x})$  be
the defining polynomial of $\cal D$, and $A=(a,b)$ where $a,b$ are
unknowns.   We consider a line $\cal L$ passing through $A$ and a
generic point $Q=(z_1,z_2)$, expressed parametrically as $L(t)= A +
t (Q-A)$. Now, we take two different points $P_1:=L(t_1),
P_2:=L(t_2)$ on $\cal L$, and we consider the algebraic set
\[ {\cal S} = \left\{ (a,b,z_1,z_2,t_1,t_2,\omega) \in \K^7 \,\left| \ \begin{array}{l} \eq_1:=g(L(t_1))=0
\\
\eq_2:=g(L(t_2))=0 \\
\eq_3:= \| L(t_1)-Q\|^2 =\|L(t_2)-Q\|^2 \\
\eq_4:= \omega \cdot (t_1-t_2)=1.
\end{array} \right\} \right. \]
$\eq_1$ and $\eq_2$ ensure that $P_1,P_2\in {\cal D}$, $\eq_3$
requires that $P_1,P_2$ are on the same circle centered at $Q$, and
$\eq_4$ guarantees that $P_1\neq P_2$. Therefore, if
$\pi:\K^7\rightarrow \K^2$ is the projection
$\pi(a,b,z_1,z_2,t_1,t_2,\omega)=(a,b)$, then the possible focuses
$A$, such that there exists $d$ for which $\frak{C}({\cal D},A,d)$
has an special component, belong to $\overline{\pi({\cal S})}$; i.e.
belong to $V(I\cap \K[a,b])$, where $I$ is the ideal generated by
$\{\eq_1,\ldots,\eq_4\}$ (see Closure Theorem in \cite{Cox} p. 122).
\begin{example}\label{example-finding-the-focus} Let $\cal D$ be the
line defined by $g(x_1,x_2)=\lambda\,x_{{1}}+\mu\,x_{{2}}+\rho$.
Then,
\[ \begin{array}{ll} \eq_1=&\lambda\, \left( a+t_{{1}} \left( z_{{1}}-a \right)  \right) +\mu\,
 \left( b+t_{{1}} \left( z_{{2}}-b \right)  \right) +\rho \\
\eq_2=&\lambda\,
 \left( a+t_{{2}} \left( z_{{1}}-a \right)  \right) +\mu\, \left( b+t_
{{2}} \left( z_{{2}}-b \right)  \right) +\rho \\ \eq_3=&\left(
a+t_{{1}}
 \left( z_{{1}}-a \right) -z_{{1}} \right) ^{2}+ \left( b+t_{{1}}
 \left( z_{{2}}-b \right) -z_{{2}} \right) ^{2}- \\ & \left( a+t_{{2}}
 \left( z_{{1}}-a \right) -z_{{1}} \right) ^{2}- \left( b+t_{{2}}
 \left( z_{{2}}-b \right) -z_{{2}} \right) ^{2} \\ \eq_4=&\omega \left( t_{{1}}-t_{
{2}} \right) -1. \end{array} \] Moreover,

\noindent
$\{\mu\,b+\lambda\,a+\rho,\lambda\,z_{{1}}+\mu\,z_{{2}}+\rho, \left(
-z_ {{2}}+b \right) ^{2} \left( t_{{1}}+t_{{2}}-2 \right) , \left(
-z_{{2} }+b \right) ^{2} \left( -2\,\omega+1+2\,\omega t_{{2}} \right)$  \\
$,\omega t_{{1}}-\omega t_{{2}} -1 \}$

\noindent is a Gr\"obner basis  of $\{\eq_1,\ldots,\eq_4\}$ w.r.t.
lexorder, with $\omega>t_1>t_2>z_1>z_2>a>b$. Hence the possible
focuses $(a,b)$ satisfy $\mu\,b+\lambda\,a+\rho=0$; i.e. they are on
the line $\cal D$. Indeed all of them are valid (see Remark
\ref{Remark-3}). \qed
\end{example}

\noindent \underline{\sf Detecting the Conchoid.} Let $\cal D$ and
$g(\bar{x})$  be as above, and let us assume that we are given know
a focus $A=(a,b)$ and we want to decide whether $\cal D$ is a
conchoid from $A$. So, we need to decide whether there exists $d\in
\K^*$ such that $\frak{C}({\cal D},A,d)$ has an special component.
For this purpose, let  $L^*(t),\eq_{1}^{*},\ldots,\eq_{4}^{*}$  be
the polynomials $L,\eq_i$  specialized at $A$. The new algebraic
set, defined by $\{\eq_{1}^{*},\ldots,\eq_{4}^{*}\}$, say ${\cal
S}^*$, is in $\K^5$. In this situation, if $\pi^*:=\K^5\rightarrow
\K^2$ where $\pi^*(z_1,z_2,t_1,t_2,\omega)=(z_1,z_2)$ and $I^*$ the
ideal generated by $\{\eq_{1}^{*},\ldots,\eq_{4}^{*}\}$, reasoning
similarly as above, the algebraic Zariski closure ${\cal
H}:=\overline{\pi^{*}({\cal S}^{*})}=V(I^*\cap \K[z_1,z_2])$
contains the special components of $\frak{C}({\cal D},A,d)$. Thus,
we may factor $\cal H$, and for each irreducible component we
compute its generic conchoid (see Remark
\ref{remark-a-la-def-conch}) to afterwards checking whether for some
$d\in \K^*$ we get $\cal D$. Note that in $\cal H$ we also may have
the curve defining the geometric locus of those points $Q$ such that
when intersecting $\cal D$ with the line passing though $A,Q$  we
get at  two points $P_1,P_2$ satisfying that $\|P_1-Q\|=\|P_2-Q\|$,
being this quantity non-constant. On the other hand, we only know
that for all, but finitely many exceptions, specializations of the
generic conchoid we get the conchoid. These exceptions can be
determined (see Example \ref{Example-1}) but the computation may be
too heavy.
\begin{example}\label{example-detecting-a-conchoid}
Let $\cal D$ be defined by
$g(\bar{x})={x_{{1}}}^{4}+2\,{x_{{2}}}^{2}{x_{{1}}}^{2}-9\,{x_{{1}}}^{2}-4\,x_{{1}
}-9\,{x_{{2}}}^{2}+12+{x_{{2}}}^{4}$ and $A=(-2,0)$. Computing a
Gr\"obner basis  of the ideal $I^*$ generated by
$\{\eq_{1}^{*},\ldots,\eq_{4}^{*}\}$  w.r.t. lexorder, with
$\omega>t_1>t_2>z_1>z_2$, we get that ${\cal H}:=V(I^*\cap
\K[z_1,z_2])$ decomposes as the union of the circle ${\cal H}_1$
defined by $ -4+{z_{{1}}}^{2}+{z_{{2}}}^{2}$ and the quartic ${\cal
H}_2$ defined by $-4-4\,z_{{1}}+
15\,{z_{{1}}}^{2}+16\,{z_{{1}}}^{3}+4\,{z_{{1}}}^{4}-{z_{{2}}}^{2}+16
\,{z_{{2}}}^{2}z_{{1}}+8\,{z_{{2}}}^{2}{z_{{1}}}^{2}+4\,{z_{{2}}}^{4}$.
The generic conchoid  $\frak{C}({\cal H}_1,A,d)$ is given by
\[ G(\bar{x},d)= -8\,{x_{{1}}}^{2}-8\,{x_{{2}}}^{2}-4\,{d}^{2}+16-4\,x_{{1}}{d}^{2}+{x_
{{1}}}^{4}+2\,{x_{{2}}}^{2}{x_{{1}}}^{2}-{x_{{1}}}^{2}{d}^{2}-{x_{{2}}
}^{2}{d}^{2}+{x_{{2}}}^{4}.\] Solving the algebraic system is $d$
provided by $g(\bar{x})=G(\bar{x},d)$ one gets that $d=\pm 1$.
Indeed ${\cal D}=\frak{C}({\cal H}_1,A,1)$. Performing the same
computations with ${\cal H}_2$ one gets that ${\cal D}$ is not a
conchoid of ${\cal H}_2$. \qed
\end{example}

\end{document}